\documentclass{amsart}

\usepackage{amscd}
\usepackage[T1]{fontenc}
\usepackage[latin1]{inputenc}
\usepackage[hypertex]{hyperref}
\usepackage[arrow,curve,matrix]{xy}
\renewcommand{\P}{\mathbb{P}}

\DeclareMathOperator{\Hom}{Hom}

\DeclareMathOperator{\Jet}{Jet}
\DeclareMathOperator{\join}{Join}
\DeclareMathOperator{\locus}{loc}
\DeclareMathOperator{\mult}{mult}
\DeclareMathOperator{\tancone}{TanCone}
\DeclareMathOperator{\Pic}{Pic}
\DeclareMathOperator{\Prolong}{Prolong}
\DeclareMathOperator{\RatCurves}{RatCurves}
\DeclareMathOperator{\Sym}{Sym}

\newcommand{\ilabel}[1]{\newcounter{#1}\setcounter{#1}{\value{enumi}}}
\newcommand{\iref}[1]{\setcounter{enumi}{\value{#1}}\labelenumi}

\theoremstyle{plain}    
\newtheorem{thm}{Theorem}[section]
\newtheorem{defn}[thm]{Definition}

\numberwithin{equation}{section}
\numberwithin{figure}{section}
\theoremstyle{plain}    
\newtheorem{cor}[thm]{Corollary}
\newtheorem{lem}[thm]{Lemma}
\theoremstyle{plain}    
\newtheorem{prop}[thm]{Proposition}
\theoremstyle{remark}
\newtheorem{fact}[thm]{Fact}
\newtheorem{claim}[thm]{Claim} 
\newtheorem{rem}[thm]{Remark}
\theoremstyle{remark}    
\newtheorem{notation}[thm]{Notation} 
\newtheorem{assumption}[thm]{Assumption}

\begin{document}

\title{Geometry of chains of minimal rational curves}

\date{\today}

\author{Jun-Muk Hwang}
\author{Stefan Kebekus}

\address{Jun-Muk Hwang, Korea Institute for Advanced Study, 207-43
  Cheongnyangni-dong, Seoul, 130-722, Korea}
\email{jmhwang@kias.re.kr}

\address{Stefan Kebekus, Mathematisches Institut, Universität zu Köln,
  Weyertal 86--90, 50931~Köln, Germany}
\email{stefan.kebekus@math.uni-koeln.de}
\urladdr{http://www.mi.uni-koeln.de/\~{}kebekus}

\thanks{J.-M.~Hwang was supported by the Korea Research Foundation
  Grant (KRF-2002-070-C00003). S.~Kebekus was supported in part by a
  Heisenberg fellowship and by the Schwerpunkt ``Globale Methoden in
  der komplexen Geometrie'' of the Deutsche Forschungsgemeinschaft.
  The work on this paper started during a visit of J.-H.~Hwang to
  Germany which was funded by the Schwerpunkt.}

\begin{abstract}
  Chains of minimal degree rational curves have been used as an
  important tool in the study of Fano manifolds. Their own geometric
  properties, however, have not been studied much.  The goal of the
  paper is to introduce an infinitesimal method to study chains of
  minimal rational curves via varieties of minimal rational tangents
  and their higher secants.  For many examples of Fano manifolds this
  method can be used to compute the minimal length of chains needed to
  join two general points.  One consequence of our computation is a
  bound on the multiplicities of divisors at a general point of the
  moduli of stable bundles of rank two on a curve.
\end{abstract}

\maketitle
\setcounter{tocdepth}{1}
\tableofcontents

\section{Introduction}

Minimal degree rational curves play an essential role in the study of
the geometry of Fano manifolds, see e.g. \cite{Hwa00},
\cite{Kebekus02a} and \cite{K96}. It is, however, frequently important
to consider not just single curves, but connected chains of them.
Among others, these were used in the proof of the boundedness of the
degrees of Fano manifolds of Picard number 1 in \cite{Nadel91}.  At
present, not much is known about the geometry of chains of curves. A
basic problem in this direction is the computation of the dimension
$d_k$ of the locus of the family of length-$k$ chains of minimal
rational curves passing through a fixed, general point of the Fano
manifold; we refer to Section~\ref{chap:spanningDim} for a precise
definition. This information could for instance be used to bound the
multiplicities of divisors at a general point of the variety, compare
Proposition~\ref{prop:3} below.  For example, our
Corollary~\ref{cor:c11} below gives a bound on the multiplicities of
divisors on the moduli space of stable rank-2 bundles on a curve.

It turns out, however, that it is not easy to compute the dimensions
$d_k$ by a direct method even for concrete examples. The goal of this
note is therefore to develop an \emph{infinitesimal} approach to this
problem.  Our main results, Theorems~\ref{thm:main1}--\ref{thm:main3},
give a lower bound on the $d_k$, which often becomes effective. In
Section~\ref{chap:examples}, we will illustrate this by computing
$d_k$ for several examples.

In addition to the intrinsic interest of computing these dimensions,
which are natural geometric invariants of Fano manifolds, our argument
shows a connection of our problem with the classical problem on the
dimensions of higher secant varieties of a projective variety,
cf.~\cite{CC02}, \cite{Zak93}. This connection is provided by the
theory of varieties of minimal rational tangents which is briefly
recalled in Section~\ref{sec:2}. In this sense, this note is another
example of the manifestation of the philosophy that the projective
geometry of varieties of minimal rational tangents governs the
geometry of the Fano variety, a phenomenon encountered repeatedly in
several works on Fano varieties.

\section{Loci of chains of rational curves}
\label{chap:spanningDim}

\subsection{Setup and definitions}

We follow the conventions of \cite{K96} and denote the normalization
of the space of irreducible, reduced rational curves on a
complex-projective variety $X$ by $\RatCurves^n(X)$.  We will
make the following assumption throughout the present paper.

\begin{assumption}\label{ass:dagger}
  Let $X$ be a uniruled complex-projective manifold and ${\mathcal K}
  \subset \RatCurves^n(X)$ be an irreducible component such that
  \begin{enumerate}
  \item \ilabel{i:covering} the members $C \in \mathcal K$ cover a
    dense subset of $X$, i.e. ${\mathcal K}$ is a dominating family of
    rational curves.
    
  \item \ilabel{i:compactness} For a general point $x \in X$, the
    subvariety
    $$
    {\mathcal K}_x = \{ C \in \mathcal K : x \in C \} \subset
    \mathcal K
    $$
    of members of ${\mathcal K}$ that pass through $x$, is compact.
    
  \item \ilabel{i:irred} For general $x \in X$, ${\mathcal K}_x$ is
    also irreducible.
  \end{enumerate}
\end{assumption}

\begin{rem}
  It is well known that~\iref{i:compactness} is satisfied if $\mathcal
  K$ is chosen such that the degrees of the member curves are minimal
  among all families of rational curves that
  satisfy~\iref{i:covering}.
  
  Assumption~\iref{i:irred}, which we pose for simplicity, is known to
  hold for essentially all examples of Fano manifolds of Picard number
  1 as long as \iref{i:covering} and \iref{i:compactness} are
  satisfied and $\dim {\mathcal K}_x \geq 1$. This is the case for
  contact Fano manifolds \cite[thm.~1.1]{Keb03} and any of the
  examples discussed in Section~\ref{chap:examples}. Examples show
  that Assumption~\iref{i:irred} need not hold when $\dim {\mathcal
    K}_x =0$. But in the latter case, the numbers corresponding to
  $d_k$ we define below are easy to compute. Thus it is reasonable to
  make the assumption~\iref{i:irred}.
  
  We refer to \cite{Hwa00} for a general discussion, and to
  \cite[thm.~5.1]{KK02} for a general irreducibility criterion.
\end{rem}

In this paper, we are principally interested in the loci of length-$k$
chains of $\mathcal K$-curves that pass through a fixed general point
$x \in X$.

\begin{defn}\label{defn:VMRT}
  For a general point $x \in X$, define inductively
  $$
  \locus^1(x) := \mbox{ closure of } \bigcup_{ C \in {\mathcal
      K}_x} C
  $$
  and
  $$
  \locus^{k+1}(x) := \mbox{ closure of } \bigcup_{\begin{array}{c}
      \mbox{general } C \in {\mathcal K}_y \\ \mbox{general } y \in
      \locus^k(x)\end{array}} C.
  $$
  Let $d_k := \dim \locus^k(x)$. We call $d_k$ the $k$-th
  \emph{spanning dimension} of the family ${\mathcal K}$.
\end{defn}

By Assumption~\ref{ass:dagger}\iref{i:irred}, $\locus^k(x)$ is an
irreducible variety.  Thus the meaning of ``general $y \in
\locus^k(x)$'' is clear in the definition of $\locus^{k+1}(x)$.  It is
clear that the spanning dimensions are independent of the choice of
the general point $x$.  Standard techniques of Mori theory give a
first estimate.

\begin{prop}\label{prop:1}
  If $C \in K$ is any curve, then $d_1 = (-K_X).C-1$ and $d_k \leq
  k((-K_X).C-1)$.
\end{prop}
\begin{proof}
  It is well-known \cite{K96} that if $x \in X$ is a general point,
  then $(-K_X).C \geq 2$ and $\dim \mathcal K_x = (-K_X).C - 2$. By
  Mori's Bend-and-Break, we have $d_1 = \dim \mathcal K_x +1$, and the
  first equality is shown.
  
  The latter inequality follows from induction. If $y \in
  \locus^{k-1}(x)$ is a general point, then
  \begin{equation*}
  \begin{split}
    \dim \locus^k(x) & \leq \dim \locus^{k-1}(x) + \dim \locus^1(y) \\
    & \leq (k-1) ((-K_X).C-1) + ((-K_X).C-1) \\
    & = k((-K_X).C-1).
  \end{split}
  \end{equation*}
\end{proof}

\subsection{Rational connectedness and the length of a manifold}

Suppose there exists an integer $\ell>0$ such that
$$
d_{\ell} = \dim X, \mbox{ but } d_{\ell -1} < \dim X. 
$$
Then we say that $X$ has {\it length} $\ell$ with respect to
${\mathcal K}$, or that $X$ is $\ell$-\emph{connected} by ${\mathcal
  K}$.  The existence of such $\ell$ is guaranteed in a number of
interesting setups.

\begin{prop}[{\cite{Nadel91}}]\label{prop:nadel}
  Under the Assumption~\ref{ass:dagger}, if $X$ is a Fano manifold of
  Picard number 1, then $X$ has length $\leq \dim X$ with respect to
  $\mathcal K$. \qed
\end{prop}

In many examples of Fano manifolds $X$ of Picard number 1, there is a
natural choice of a family ${\mathcal K}$ that satisfies the
Assumption~\ref{ass:dagger}. The length $\ell$ and the
spanning-dimensions $d_1, \ldots, d_{\ell}$ can then be regarded as
natural numerical invariants of $X$. In addition to its intrinsic
interest, these numbers are useful in understanding the geometry of
$X$. One use of the length is the following generalization of Nadel's
product theorem \cite{Nadel91}.

\begin{prop}\label{prop:3}
  Under the Assumption~\ref{ass:dagger}, let $V$ be a vector bundle on
  $X$. Consider a general curve $C \in \mathcal K$, let $\nu: {\mathbb
    P}_1 \to C$ be its normalization, and write
  \begin{equation}
    \label{eq:splittypeofv}
    \nu^* (V) \cong {\mathcal O}_{\mathbb P_1}(a_1) \oplus \cdots
    \oplus {\mathcal O}_{\mathbb P_1}(a_r)
  \end{equation}
  for some integers $a_1 \geq \cdots \geq a_r$.
  
  If $X$ is of length $\ell$ with respect to $\mathcal K$, $x \in X$
  is a general point and $\sigma \in H^0(X, V)$ is any non-zero
  section, then the order of vanishing of $\sigma$ at $x$ satisfies
  $$
  \mult_x(\sigma) \leq \ell \cdot a_1.
  $$
\end{prop}

\begin{rem}
  The proof of Proposition~\ref{prop:3} employs higher-order jet
  bundles and the associated short exact sequences. The reader who is
  not familiar with those might wish to consult \cite[sect.~4.1 and
  4.6]{KS72}. Basic fact about standard rational curves can be found
  in Section~\ref{sect:stdcurves} below.
\end{rem}

\begin{proof}[Proof of Proposition~\ref{prop:3}]
  Set $x_0 := x$ and let $x_{i+1}$ be a general point of
  $\locus^1(x_i)$.  We may assume that there exists a general member
  $C_i$ of ${\mathcal K}$ containing both $x_i$ and $x_{i+1}$. Let
  $m_i := \mult_{x_i} (\sigma)$ and let $\nu_i : \mathbb P_1 \to C_i$
  be the normalization morphisms. Note that $m_{\ell} =0$ by the
  definition of the length $\ell$.
  
  Since $x_i$ was generically chosen, the section $\sigma$ vanishes
  along $C_i$ to order exactly $m_{i+1}$, and $m_{i+1} \leq m_i$.
  This has two consequences. First, if
  $$
  \sigma_{m_{i+1}} := \Prolong^{m_{i+1}}(\sigma) \in H^0(X,
  \Jet^{m_{i+1}}(V))
  $$
  is the $m_{i+1}$th jet of $\sigma$, then $\sigma_{m_{i+1}}$ has a
  zero of order $m_i - m_{i+1}$ at the point $x_i$. Second, if
  $$
  \begin{CD}
    0 @>>> V \otimes \Sym^{m_{i+1}}\Omega^1_X @>>>
    \Jet^{m_{i+1}}(V) @>{\alpha}>> \Jet^{m_{i+1}-1}(V) @>>> 0
  \end{CD}
  $$
  is the $m_{i+1}$th jet-sequence of $V$, then
  $\nu_i^*(\alpha(\sigma_{m_{i+1}})) = 0$, and $\sigma_{m_{i+1}}$
  therefore defines a non-zero section
  $$
  \tau \in H^0 \left(\mathbb P_1, \nu_i^*(V \otimes
    \Sym^{m_{i+1}}\Omega^1_X)\right)
  $$
  which has a zero of order $m_i - m_{i+1}$ at the point $x_i$.  We
  conclude with an analysis of the splitting type of $\nu_i^*(V
  \otimes \Sym^{m_{i+1}}\Omega^1_X)$. To this end, write
  $$
  \nu_i^*(V \otimes \Sym^{m_{i+1}}\Omega^1_X) \cong {\mathcal
    O}_{\mathbb P_1}(b_1) \oplus \cdots \oplus {\mathcal O}_{\mathbb
    P_1}(b_k)
  $$
  for integers $b_1 \geq \cdots \geq b_k$. By
  Fact~\ref{fact:mostcurvesstd}, the general curves $C_i$ are
  standard, and
  $$
  \nu_i^* (\Omega^1_X) \cong {\mathcal O}_{\mathbb P_1}(-2) \oplus
  {\mathcal O}_{\mathbb P_1}(-1)^{\oplus a} \oplus {\mathcal
    O}_{\mathbb P_1}^{\oplus b}.
  $$
  This, together with the isomorphism \eqref{eq:splittypeofv}
  implies that $b_1 \leq a_1$, and the existence of the section $\tau$
  implies that $m_i - m_{i+1} \leq b_1 \leq a_1$ for all $i$. Thus
  $$
  \mult_x(\sigma) = m_0 = \sum_{i=0}^{\ell-1}(m_i-m_{i+1}) \leq
  \ell \cdot a_1.
  $$
\end{proof}

\section{Statement of the main result}
\label{sec:2}

The main results of this paper relate the dimensions $d_k$ to the
projective geometry of the variety of minimal rational tangents. For
the reader's convenience, we briefly recall the relevant definitions
first.

\subsection{Varieties of minimal rational tangents and higher secants}

\begin{fact}
  Under the Assumption~\ref{ass:dagger}, if $x \in X$ is a general
  point, then the normalization $\hat{\mathcal K}_x$ of ${\mathcal
    K}_x$ is smooth. It is shown in \cite[thm.~3.1]{Kebekus02a} that
  the rational map
  $$
  \begin{array}{ccl}
    {\mathcal K}_x & \dasharrow & \P(T_X|_x^\vee) \\
    C & \mapsto & \P(T_C|_x^\vee)
  \end{array}
  $$
  that sends a curve through $x$ to its tangent direction at $x$,
  extends to a finite morphism $\tau_x: \hat{\mathcal K}_x \to
  \P(T_X|_x^\vee)$.
\end{fact}

\begin{defn}\label{def:Cx}
  For general $x \in X$, let ${\mathcal C}_x$ be the image
  $\tau_x(\hat{\mathcal K}_x)$ in $\P(T_X|_x^\vee)$ and call it the
  variety of minimal rational tangents at $x$.
\end{defn}

The variety of minimal rational tangents comes with a natural
embedding into $\P(T_X|_x^\vee)$. Its invariants as a projective
variety can therefore be used to define natural invariants of the
manifold $X$ and the family $\mathcal K$. To this end, let us recall
the notion of higher secant varieties.

\begin{defn}[{\cite{CC02} and \cite{Zak93}}]\label{def:highersecant}
  Let $Y \subset \P_m$ be an irreducible projective subvariety of
  dimension $p$. The $k$-th {\it secant variety} $S^k Y$ of $Y$ is
  defined to be the closure of the union of $k$-dimensional linear
  subspaces of $\P_m$ determined by general $k+1$ points on $Y$.
\end{defn}

\begin{fact}
  In the setup of Definition~\ref{def:highersecant}, let $m'$ be the
  dimension of the linear subspace spanned by $Y$. Then the dimension
  of $S^k Y$ satisfies
  $$
  \dim(S^k Y) \leq \min\{ m', (k+1)(p+1)-1 \}.
  $$
\end{fact}

\begin{notation}
  We say that $Y$ is $k$-\emph{defective} when the above inequality
  is strict.
\end{notation}

\begin{notation}\label{not:VMRT}
  As a subvariety of the projective space $\P(T_X|_x^\vee)$, we can
  consider the $k$-th secant variety $S^k{\mathcal C}_x$ of ${\mathcal
    C}_x$. For $k \geq 1$, let ${\mathcal C}_x^{k} \subset
  \P(T_X|_x^\vee)$ be the tangent cone of $\locus^{k}(x)$ at $x$ ---we
  refer to Section~\ref{sec:tancone} for a discussion of tangent
  cones.
  
  Finally, let $\mathcal C$ and $S^k \mathcal C \subset \P(T_X^\vee)$
  be the closure of the union of the $\mathcal C_x$ and $S^k \mathcal
  C_x$ for all general points $x \in X$. We call $\mathcal C$ the
  variety of minimal rational tangents.
\end{notation}

\begin{rem}\label{rem:closednessOfC}
  It is clear that $d_k = \dim {\mathcal C}_x^{k} +1$ for general
  points $x \in X$.  If $x \in X$ is a general point, it is likewise
  obvious that $\mathcal C \cap \P(T_X|_x^\vee) = \mathcal C_x$ and
  $S^k \mathcal C \cap \P(T_X|_x^\vee) = S^k\mathcal C_x$.
\end{rem}

\begin{defn}[{\cite[I (1.1.1)]{Zak93}}]
  Let $A,B \subset \P_N$ be two irreducible subvarieties of the
  projective space.  Define $\join(A,B)$ to be the closure of the
  union of lines joining a general point of $A$ to a general point of
  $B$.
\end{defn}

\begin{rem}[{\cite[V (1.2.5)]{Zak93}}]
  In the notation of the definition above, $S^1 A = \join(A,A)$ and
  $S^kA = \join(A, S^{k-1}A)$.
\end{rem}

\begin{notation}
  Let $\ell \subset X$ be a rational curve and $\eta: \mathbb P_1 \to
  \ell$ its normalization. We say that $\ell$ has \emph{nodal
    singularities} at a point $x \in \ell$, if $\eta^{-1}(x)$ contains
  more than one point. We say that $\ell$ has \emph{immersed
    singularities} at $x$, if the morphism $\eta$ has rank one at all
  points of the preimage $\eta^{-1}(x)$.
\end{notation}

\subsection{Main Results}

The following are the main results of this paper. In each of the
statements, we maintain the Assumptions~\ref{ass:dagger} and let $x
\in X$ be a general point.

\begin{thm}\label{thm:main1}
  We have ${\mathcal C}_x^1 \subset S^1 \mathcal C_x$. If none of the
  curves associated with $\mathcal K_x$ has a nodal singularity at
  $x$, then ${\mathcal C}_x^1 = {\mathcal C}_x$.
\end{thm}

\begin{thm}\label{thm:main2}
  For each $k\geq 1$, we have $\join(\mathcal C_x, \mathcal C_x^k)
  \subseteq \mathcal C_x^{k+1}$.  In particular, $S^k {\mathcal C}_x
  \subseteq {\mathcal C}_x^{k+1}$ and $d_{k+1} \geq \dim S^k{\mathcal
    C}_x + 1$.
\end{thm}

As an immediate corollary to Theorem~\ref{thm:main2} and
Proposition~\ref{prop:1}, we obtain the following.

\begin{cor}
  If ${\mathcal C}_x$ is not $k$-defective, then $d_{k+1} = \dim
  S^k{\mathcal C_x}+1$. \qed
\end{cor}

\begin{thm}\label{thm:main3}
  If there exists an embedding $X \subset \mathbb P_N$ such that
  curves associated with $\mathcal K$ are mapped to lines in $\mathbb
  P_N$, then $S^1 {\mathcal C}_x = {\mathcal C}_x^2$. In particular,
  $d_2 = \dim S^1 {\mathcal C}_x +1$.
\end{thm}

\begin{rem}
  In Section~\ref{sect:contactmfs} we will see examples where $S^k
  {\mathcal C}_x \neq {\mathcal C}_x^{k+1}$ for $k \geq 2$, even when
  $X$ is embedded into $\P_N$ in a way such that curves associated
  with $\mathcal K$ are lines in $\P_N$.
\end{rem}

\begin{cor}
  Suppose $X$ is a Fano manifold of Picard number 1 and $\dim X \geq
  3$ such that the ample generator $H$ of $\Pic(X)$ is very ample and
  $-K_X = i H$ for some
  \begin{equation}
    \label{eq:bigindex}
    i > \frac{2}{3} \dim X + \frac{1}{3}.    
  \end{equation}
  Let ${\mathcal K}$ be the set of lines lying on $X$ under the
  embedding $X \subset \P(H^0(X, H))$.  Then $X$ is 2-connected by
  ${\mathcal K}$. A hyperplane section through a general point $x \in
  X$ has multiplicity $\leq 2$ at $x$, or equivalently, the second
  fundamental form of $X$ at $x$ in $\P(H^0(X, H))$ is surjective to
  the normal space.
\end{cor}

\begin{proof}
  By \cite[p.~348]{Hwa00}, ${\mathcal C}_x$ is a smooth subvariety of
  $\P(T_X|_x^\vee)$ of dimension $i-2$. It follows that ${\mathcal
    C}_x$ is irreducible, and \cite[thm.~2.5]{Hwa00} implies that it
  is nondegenerate. By \cite[II.2.1]{Zak93}, the inequality
  \eqref{eq:bigindex} then yields that the secant variety is the full
  ambient space, $S^1{\mathcal C}_x = \P(T_X|_x^\vee)$.
  Theorem~\ref{thm:main3} therefore implies that $X$ is 2-connected by
  lines. The last sentence follows from Proposition~\ref{prop:3}.
\end{proof}

\section{Facts used in the proof of Theorems~\ref{thm:main1}--\ref{thm:main3}}

The proof of the main results relies on a number of facts that are
scattered throughout the literature. For the reader's convenience, we
have gathered those here.

\subsection{Tangent cones}\label{sec:tancone}

Throughout the proof of Theorems~\ref{thm:main1}--\ref{thm:main3}, we
use the following two descriptions of the projectivized tangent cone
$\tancone_{Y,x} \subset \P(T_X|_x^\vee)$ of a subvariety $Y \subset X$
at a point $x \in Y$.  For an elementary discussion of holomorphic
arcs, of tangent cones and for a proof of these descriptions, see
\cite[lect.~20]{Harris95}.

\begin{enumerate}
\item As a set, the tangent cone is the union of the tangent lines to
  holomorphic arcs $\gamma: \Delta \to Y$ that are centered about $x$.
  Recall that if the arc $\gamma$ is not smooth at $0$, then the
  tangent line at 0 is just the limit of the tangent lines to smooth
  locus. For example, in Lemma~\ref{lem:tanofopen}, when we say that
  $\gamma$ has tangent $\alpha$ we do not mean that $\gamma$ is smooth
  there.
  
\item \ilabel{lbl:blowup} Consider the blow up $\hat X \to X$ of the
  point $x$, and identify the exceptional divisor $E \subset \hat X$
  with $\mathbb P(T_X|_x^\vee)$. If $\hat Y \subset Y$ is the strict
  transform of $Y$, then the projectivized tangent cone of $Y$ at $x$
  is exactly the intersection $\hat Y \cap E$.
\end{enumerate}

The following lemma is immediate from \iref{lbl:blowup}.

\begin{lem}\label{lem:tanofopen}
  Let $Y_o \subset Y$ be a Zariski open subset of an irreducible
  subvariety $Y \subset X$. For any point $y \in Y$ and a general
  point $\alpha$ of any component $T$ of $\tancone_{Y,y}$, there
  exists an arc $\gamma$ on $Y$ with tangent $\alpha$ at $y$ such that
  $\gamma \cap Y_o \neq \emptyset$. \qed
\end{lem}

\begin{cor}\label{cor:tancone}
  Let $Y \subset X$ be an irreducible subvariety, and set 
  $$
  \mathcal T := \bigcup_{y \in Y} \tancone_{Y,y} \subset \P(T_X^\vee).
  $$
  Then the reduced constructible set $\mathcal T$ is also
  irreducible.
\end{cor}

\begin{proof}
  Suppose $\mathcal T$ is not irreducible. Since for a general point
  $y \in Y$, the tangent cone $\tancone_{Y,y}$ is a linear subspace in
  $\P(T_X|_y^\vee)$, then there exists a component $\mathcal Z \subset
  \mathcal T$ whose image in $Y$ under the natural projection is a
  proper subvariety.  Let $Y_o$ be its complement.  Choose a general
  point $\alpha \in \mathcal Z$ which does not lie in other components
  of $\mathcal T$. By Lemma~\ref{lem:tanofopen}, there exists an arc
  $\gamma$ tangent to $\alpha$ which intersects $Y_o$.  The lift of
  $\gamma$ to $\P(T_X^\vee)$ lies in $\mathcal T$ and contains
  $\alpha$. Consequence: the lift lies in $\mathcal Z$, a
  contradiction to the choice of $Y_o$.
\end{proof}

\subsection{Vector bundles, free and standard rational curves}
\label{sect:stdcurves}

\begin{defn}\label{def:standard}
  A rational curve $\ell \in \mathcal K$ is ``free'' if the
  pull-back of the tangent bundle $T_X$ via the normalization $f: \P_1
  \to \ell$ is isomorphic to
  $$
  f^*(T_X) \cong \bigoplus_{i=1 \ldots \dim X} {\mathcal
    O}_{\P_1}(a_i)
  $$
  where all $a_i \geq 0$.  The curve $\ell$ is ``standard'' if
  there are numbers $a$, $b$ such that
  $$
  f^*(T_X) \cong {\mathcal O}_{\P_1}(2) \oplus {\mathcal
    O}_{\P_1}(1)^{\oplus a} \oplus {\mathcal O}_{\P_1}^{\oplus b}.
  $$
\end{defn}

\begin{fact}[{\cite[sect.~1.2]{Hwa00}}]\label{fact:mostcurvesstd}
  The sets of free and standard curves in $\mathcal K$ are open.  If
  $x \in X$ is a general point and $\ell \in \mathcal K$ corresponds
  to a curve that contains $x$, then the curve associated with $\ell$
  is free.  If $\ell$ corresponds to a general point of $\mathcal
  K_x$, the curve is standard. \qed
\end{fact}

\begin{lem}\label{lem:Ksmooth}
  If $\ell \in \mathcal K$ corresponds to a free curve, then the
  parameter space $\mathcal K$ is smooth at $\ell$, and $\ell$ is a
  fiber of the natural morphism
  $$
  \mathcal K \to \text{Chow-variety of curves in $X$.}
  $$
  In particular, since freeness is an open property, the morphism
  is generically injective.
\end{lem}
\begin{proof}
  Both assertions follow from \cite[thm.~II.2.15]{K96} if we show that
  the $\Hom$-scheme $\Hom_{bir}(\P_1, X)$ is smooth in the
  neighborhood of a morphism $f : \P_1 \to X$ that parameterizes the
  curve associated with $\ell$. But since that curve is free, the
  $\Hom$-scheme is indeed smooth at $f$, \cite[cor.~II.3.5.4]{K96}.
\end{proof}

\begin{notation}
  We have already used the fact that any vector bundle $E$ on $\mathbb
  P_1$ decomposes as a sum of line bundles, $E = \bigoplus \mathcal
  O_{\mathbb P_1}(a_i)$. While this decomposition is not unique in
  general, the positive sub-bundle
  $$
  E^+ := \bigoplus_{\{i\,|\, a_i > 0\}} {\mathcal O}_{\P_1}(a_i)
  \subset E
  $$
  is in fact independent of the choice of the decomposition.
  Throughout this paper, we attach the index ``$+$'' to the names of
  vector bundles on $\mathbb P_1$ to denote their positive parts.
\end{notation}

The following Lemma is immediate.

\begin{lem}\label{lem:splitOfStd}
  Let $\ell$ be a smooth standard rational curve on $X$ and $\vec v$ a
  vector in the positive part of $T_X|_\ell^+$, not tangent to $\ell$.
  Then there exists a unique subbundle $F$ of rank 2 in $T_X|_\ell$
  which contains $\vec v$ and splits as ${\mathcal O}_{\P_1}(2) \oplus
  {\mathcal O}_{\P_1}(1)$.
  
  When $\ell$ is a line in $\P_N$, then $F$ is the restriction of the
  tangent bundle of the plane determined by $\ell$ and $\vec v$.  \qed
\end{lem}

\begin{fact}[{\cite[prop.~II.3.4]{K96}}]\label{faxt:T+andS1}
  Let $x \in X$ be a general point, and $\ell \in \mathcal K$
  correspond to a standard rational curve that contains $x$. Then the
  positive part of the restriction $T_X|_{\ell}$ is contained in the
  first secant variety $S^1\mathcal C_x$. More precisely, if $\ell$ is
  parameterized by a morphism $f: \mathbb P_1 \to \ell$ and if $y \in
  f^{-1}(x)$ is any point, then $f$ is an immersion and
  $$
  \mathbb P\left( f^*(T_X)^+|_y^\vee\right) \subset f^*(S^1
  \mathcal C_x) \subset \mathbb P\left( f^*(T_X)|_y^\vee\right).
  $$
\end{fact}

\subsection{Varieties that are covered by lines}

If there exists an embedding $X \subset \P_N$ such that curves
associated with $\mathcal K$ correspond to lines in $\P_N$, then
$\mathcal K$ is compact and all $\mathcal K$-curves are smooth. Many
of the statements mentioned previously can be improved.

\begin{lem}[{\cite[proof of prop.~1.5]{Hwa00}}]\label{lem:almostallcurvesstd}
  If there exists an embedding $X \subset \P_N$ such that $\mathcal
  K$-curves are lines, then free lines are standard. In particular, by
  Fact~\ref{fact:mostcurvesstd}, if $x \in X$ is a general point and
  if $\ell \subset X$ is a line which contains $x$, then $\ell$ is
  standard.
\end{lem}
\begin{proof}
  Let $\ell \subset X$ be a free line. We may write
  $$
  T_X|_{\ell} \cong {\mathcal O}_{\P_1}(a_1) \oplus \cdots \oplus
  {\mathcal O}_{\P_1}(a_n)
  $$
  with $a_1 \geq \cdots \geq a_n \geq 0$. We know that $a_1 \geq 2$
  because the smoothness of $\ell$ implies that $T_\ell$ injects into
  $T_X|_\ell$.  On the other hand, $T_X|_\ell$ is a subbundle of
  $$
  T_{\P_N}|_\ell \cong {\mathcal O}_\ell(2) \oplus [{\mathcal
    O}_\ell(1)]^{\oplus N-1}.
  $$
  Thus $a_1=2$ and $1 \geq a_2 \geq \cdots \geq a_n \geq 0$ and $\ell$
  is standard.
\end{proof}

\begin{prop}\label{prop:smx}
  Let $\mathcal C \subset \P(T_X^\vee)$ be the variety of minimal
  rational tangents that was introduced in Notation~\ref{not:VMRT} and
  let $p : \mathcal C \to X$ be the canonical projection.  If there
  exists an embedding $X \subset \P_N$ such that $\mathcal K$-curves
  are lines, if $\ell \in \mathcal K$ corresponds to a free line and
  if $x \in X$ is any point on that line, then $\mathcal C|_x$ is
  smooth at $\vec v := \mathbb P(T_\ell |_x^\vee)$ and the projective
  tangent space to $\mathcal C|_x$ at $\vec v$ is exactly
  $\P\left(T_X|_\ell^+|_x^\vee\right)$.
\end{prop}
\begin{proof}
  Let $U$ be the universal $\P_1$-bundle over $\mathcal K$, and let
  $\iota : U \to X$ be the universal morphism. Since all curves
  involved are smooth, the restriction of the tangent morphism
  $T\iota$ to the relative tangent bundle $T_{U|\mathcal K}$ yields a
  morphism $\tau$ which factors the universal morphism as follows.
  $$
  \xymatrix{ & & {\mathbb P(T_X^\vee)} \ar[d]^{p} \\ {U} \ar[rr]^{\iota}
    \ar[d]^{\pi}_{\P_1 \text{-bundle}} \ar@/^/[rru]^{\tau} & & {X}
    \\
    {{\mathcal K}}}
  $$
  The image of the morphism $\tau$ is then exactly $\mathcal C$.
  Since lines are uniquely determined by their tangent vectors,
  Lemma~\ref{lem:Ksmooth} implies that $\tau^{-1}(\vec v)$ is a single
  smooth point, say, $y \in U$. The claim that $\mathcal C|_x$ is
  smooth at $\vec v$ thus follows if we show that both $\iota$ and
  $\tau$ have maximal rank at $y$. The morphism $\iota$ has maximal
  rank because $\ell = \pi(y)$ corresponds to a free curve ---see
  \cite[cor.~II.3.5.4 and thm.~II.2.15]{K96}. The fact that $\tau$ has
  maximal rank follows from Lemma~\ref{lem:almostallcurvesstd}, which
  asserts that the curve of $\ell$ is standard, and from the
  description of the tangent map of the universal morphism $\iota$,
  \cite[thm.~II.2.15 and prop.~II.3.4]{K96}.
  
  The latter reference also explains why the tangent space to
  $\mathcal C|_x$ at $\vec v$ is $\P\left(T_X|_\ell^+|_x^\vee\right)$.
\end{proof}

\section{Proof of Theorems~\ref{thm:main1}--\ref{thm:main3}}
\label{sec:pomainth}

\subsection{Proof of Theorem~\ref{thm:main1}}

For clarity, we show the two statements of Theorem~\ref{thm:main1}
separately.

\subsubsection{Proof of the inclusion ${\mathcal C}_x^1 \subset
  S^1 \mathcal C_x$}

To shorten the notation, set $Y := \locus^1(x)$. Observe that it
follows from Assumption~\ref{ass:dagger}\iref{i:irred} that $Y$ is
irreducible.

Since $x$ is chosen to be a general point of $X$, we have the equality
$S^1\mathcal C_x = S^1 \mathcal C|_x$. In view of
Corollary~\ref{cor:tancone}, to show that $\tancone_{Y,x} \subset S^1
\mathcal C_x$, it suffices to show that for a general (hence smooth)
point $y \in Y$, the inclusion $\P(T_Y|_y^\vee) \subset S^1 \mathcal
C$ holds.

To this end, consider the diagram 
$$
\xymatrix{ {U_x} \ar[r]^(.5){ \iota_x} \ar[d]^{
    \pi_x}_{\txt{\scriptsize $\mathbb P_1$-bundle}} & {X} \\
{ {\mathcal K}_x} }
$$
where $\pi_x$ is the restriction of the universal $\P_1$-bundle,
and $\iota_x$ is the universal morphism which surjects onto
$\locus^1(x)$. Let $\hat y \in U_x$ be a general point of $U_x$.

If $\hat \ell := \pi_x^{-1}\pi_x(\hat y)$ is the fiber through $\hat
y$, then $\iota_x(\hat \ell)$ is a standard rational curve by
Fact~\ref{fact:mostcurvesstd}. In this setup generic smoothness of
$\iota_x$ and \cite[II.3.4]{K96} assert that
\begin{enumerate}
\item the morphism $\iota_x$ is smooth at $\hat y$, and
  
\item the tangent map $T\iota_x : T_{U_x}|_{\hat y} \to
  \iota_x^*(T_X)|_{\hat y}$ has its image in $(\iota_x^*(T_X)|_{\hat
    \ell})^+$.
\end{enumerate}
The inclusion 
$$
\P\left(T_Y|_{\iota_x(\hat y)}^\vee\right) \subset S^1\mathcal
C_{\iota_x(\hat y)} \subset S^1\mathcal C
$$
then follows from Fact~\ref{faxt:T+andS1}. \qed

\subsubsection{Proof of the equality ${\mathcal C}_x^1 = {\mathcal C}_x$}
\label{sec:xy}

Let $\hat {\mathcal K}_x$ be the normalization of $\mathcal K_x$, let
$\pi_x : U_x \to \hat {\mathcal K}_x$ be the pull-back of the
universal $\mathbb P_1$-bundle, and let $\iota_x : U_x \to X$ be the
universal morphism. The blow up $\hat X \to X$ of the point $x$ then
yields a diagram
$$
\xymatrix{ & & {\hat X} \ar[d]^{\txt{\scriptsize blow-up of $x$}}
 \\ { U_x} \ar[rr]^{ \iota_x} \ar[d]^{
 \pi_x}_{\txt{\scriptsize $\mathbb P_1$-bundle}} \ar@{-->}@/^/[rru]^{\hat
 \iota_x} & & {X}
 \\
 {\hat {\mathcal K}_x} }
$$
The assumption that none of the curves associated with $\mathcal
K_x$ has a nodal singularity implies that the scheme-theoretic
preimage $\sigma = \iota_x^{-1}(x) \subset U_x$ is supported on a
section over $\hat {\mathcal K}_x$. Recall
from~\cite[thm.~3.3]{Kebekus02a} that the curves associated with
$\mathcal K_x$ are either smooth or have immersed singularities at
$x$. Consequence: the preimage $\sigma$ is reduced, and it is a
Cartier-divisor in $U_x$. The universal property of the blow-up,
\cite[prop.~II.7.14]{Ha77}, then states that the rational map $\hat
\iota_x$ is actually a morphism. The composition
$$
\hat \iota_x \circ (\pi_x|_\sigma)^{-1} : \hat {\mathcal K}_x \to
E\cong \mathbb P(T_X|_x^\vee)
$$
clearly equals the tangent morphism $\tau_x$, and the equality
${\mathcal C}_x^1 = {\mathcal C}_x$ follows from the description of
the tangent cone that we gave in section~\ref{sec:tancone} above. \qed

\subsection{Proof of Theorem~\ref{thm:main2}}
\label{sec:proof_of_THM_ii}

We will prove Theorem~\ref{thm:main2} by induction.

\subsubsection{Start of induction, $k=0$}

If $k=0$, the statement of Theorem~\ref{thm:main2} reduces to
$\mathcal C_x \subseteq \mathcal C_x^1$. That, however, follows by
definition.

\subsubsection{Inductive step}

Our aim here is to show that $\join ({\mathcal C}_x, {\mathcal C}_x^k)
\subset {\mathcal C}_x^{k+1}$. We assume that the inclusion $\join
({\mathcal C}_x, {\mathcal C}_x^{k-1}) \subset {\mathcal C}_x^{k}$ is
already established.

Let $\vec v, \vec w \in T_X|_x \setminus \{0\}$ be two arbitrary
linearly independent tangent vectors whose classes $[\vec v], [\vec w]
\in \mathbb P(T_X|_x^\vee)$ are in $\mathcal C^k_x$ and $\mathcal
C_x$, respectively. In order to prove Theorem~\ref{thm:main2}, we need
to show that the line $\ell \in \mathbb P(T_X|_x^\vee)$ through $[\vec
v]$ and $[\vec w]$ is contained in $\mathcal C_x^{k+1}$, the tangent
cone to $\locus^{k+1}(x)$.

It follows from the assumptions that we can find an arc $\gamma:
\Delta \to \locus^k(x)$, where $\Delta$ is a unit disk, whose tangent
line at $0\in \Delta$ is spanned by $\vec v$.  Similarly to the
argumentation in section~\ref{sec:xy} above, let $\pi:U \to \mathcal
K$ be the universal $\mathbb P_1$ bundle and $\iota: U\to X$ the
evaluation morphism. Since $[\vec w] \in \mathcal C_x$, we find a
closed point $y \in \iota^{-1}(x)$ such that $\pi(y)$ corresponds to a
curve with tangent line $[\vec w]$, i.e.  such that
$$
\P\left( T\iota \left(T_{U/\mathcal K}|_y \right)^\vee \right) =
[\vec w] \in \mathbb P\left( T_X|_x^\vee \right).
$$
Recall from \cite[cor.~II.3.5.4]{K96} that the evaluation morphism
$\iota$ is smooth in a neighborhood of $y$. Consequence: after
shrinking $\Delta$, if necessary, we can find an arc $\tilde \gamma:
\Delta \to U$ such that
\begin{enumerate}
\item $\tilde \gamma$ is a lifting of $\gamma$, i.e.~$\iota \circ
  \tilde \gamma = \gamma$, 

\item $\tilde \gamma(0)=y$, and 
  
\item \ilabel{tgamma:finite} for general $s \in \Delta$, the curve
  associated with $(\pi \circ \tilde \gamma)(s) \in \mathcal K$ does
  not contain the point $x$, i.e.~$(\pi \circ \tilde \gamma)(\Delta)
  \not \subset \mathcal K_x$.
\end{enumerate}
Base change and blowing up will yield a diagram as follows.
$$
\xymatrix{
  {\hat U_\gamma} \ar[rrrr]^{\hat \iota_\gamma} 
  \ar[d]_{\text{blow-up $\beta_\gamma$ of $\iota_\gamma^{-1}(x)$}} &&&&
  \hat X \ar[d]^{\text{blow-up $\beta$ of $x$}} \\
  {U \times_{\mathcal K} \Delta} \ar[rr] 
  \ar[d]^{\pi_{\gamma}}
  \ar@/^0.5cm/[rrrr]^{\iota_\gamma} & &
  {U} \ar[rr]_{\iota} \ar[d]_{\pi} & &
  {X}  \\
  {\Delta} \ar[rr]_{\pi \circ \tilde \gamma} \ar@/^0.3cm/[u]^{\sigma}
  & & {\mathcal K} 
}
$$
Here $\pi$ and $\pi_\gamma$ are $\mathbb P_1$-bundles, and $\sigma$
is a section that satisfies $\iota_\gamma \circ \sigma = \gamma$.

\begin{notation}
  For brevity, set $U_\gamma := U \times_{\mathcal K} \Delta$.
  Property~\iref{tgamma:finite} of $\tilde \gamma$ ensures that all
  components of the preimage $\iota_\gamma^{-1}(x)$ are
  zero-dimensional. Let $y_\gamma \subset \iota_\gamma^{-1}(x)$ be the
  component of the scheme-theoretic preimage of $x$ that is supported
  at the point $\sigma(0)$.
  
  Let $E \subset \hat X$ be the exceptional divisor of $\beta$, and
  $E_\gamma \subset \hat U_\gamma$ be the unique irreducible and
  reduced component of the exceptional divisor of $\beta_\gamma$ that
  lies over $\sigma(0)$.
\end{notation}

\begin{rem}\label{rem:loc}
  By Lemma~\ref{lem:tanofopen}, we may assume that $\gamma$ is an arc
  in $\locus^k(x)$ such that $\iota_\gamma(U_\gamma) \subset
  \locus^{k+1}(x)$.
\end{rem}

Remark~\ref{rem:loc} and the description of the tangent cone in terms
of blowing-up together guarantee that the image $\hat
\iota_\gamma(E_\gamma)$ is contained in the tangent cone $\mathcal
C_x^{k+1}$ to $\locus^{k+1}(x)$. Since $\hat \iota_\gamma(E_\gamma)$
contains both $[\vec w]$ and $[\vec v]$, the inclusion $\ell \subset
\mathcal C_x^{k+1}$ is proved as soon as we show that $\hat
\iota_\gamma(E_\gamma)$ is a line in $E \cong \mathbb P(T_X|_x^\vee)$.
Since the normal bundle of $E$ in $\hat X$ is $N_{E/\hat X} \cong
\mathcal O_{\mathbb P_{n-1}}(-1)$, it suffices to show the following
claim.

\begin{claim}
  The intersection number $\hat \iota_\gamma^*(E).E_\gamma$ is $-1$.
\end{claim}
\begin{proof}
  Recall from~\cite[thm.~3.3]{Kebekus02a} that the curves associated
  with $\mathcal K_x$ have at worst immersed singularities at $x$. One
  consequence of this is that the scheme-theoretic intersection of the
  (possibly non-reduced) point $y_\gamma$ and the (reduced) curve
  $\pi_\gamma^{-1}(0)$ is a reduced point. We can therefore find a
  positive integer $n \in \mathbb N$ and a bundle coordinate $t$ on an
  open neighborhood $\Omega = \Omega(y_\gamma) \subset U_\gamma$ such
  that $y_\gamma$ can be expressed in terms of $t$ and the coordinate
  $s$ on $\Delta$ as follows
  $$
  y_\gamma = \{ (t,s) \in \Omega \, | \, t=0, s^n=0 \}.
  $$
  
  The blow-up of this ideal can be written down explicitly. An
  elementary computation reveals that the blow-up $\hat U_\gamma$ has
  a $\mathbb Q$-factorial singularity of type $A_{n-1}$, and that the
  following two equations of $\mathbb Q$-Cartier divisors and their
  intersection numbers hold. First:
  \begin{align*}
    \hat \iota_\gamma^*(E) & = \beta_\gamma^{-1} (y_\gamma) +
    \text{other
      components that do not intersect } E_\gamma \\
    & = n \cdot E_\gamma + \text{other components that do not
      intersect } E_\gamma
  \end{align*}
  Second:
  $$
  (n \cdot E_\gamma) \text{ . strict transform of } \pi_\gamma^{-1}
  (0) = 1.
  $$
  Since $(n \cdot E_\gamma).(\pi_\gamma \circ \beta_\gamma)^{-1}(0) =
  0$ and 
  \begin{align*}
    (\pi_\gamma \circ \beta_\gamma)^{-1}(0) & = E_\gamma +
    \text{strict transform of } \pi_\gamma^{-1} (0)\\
    & \quad + \text{other components that do not intersect } E_\gamma,
  \end{align*}
  the claim then follows immediately.
\end{proof}

This ends the proof of Theorem~\ref{thm:main2}. \qed

\subsection{Proof of Theorem~\ref{thm:main3}}

The inclusion $S^1 \mathcal C_x \subseteq \mathcal C_x^2$ has been
shown in section~\ref{sec:proof_of_THM_ii} above, and it remains to
show the opposite inclusion $\mathcal C_x^2 \subseteq S^1 \mathcal
C_x$. As the proof is somewhat lengthy, we subdivide it into steps.
Throughout the proof we will fix an embedding $X \subset \P_N$ such
that $\mathcal K$-curves become lines.

\subsubsection*{Step 1: General Setup}

Let $\vec v \in \mathcal C^2_x \subset \P(T_X|_x^\vee)$ be a general
element of an irreducible component of $\mathcal C^2_x$. We will need
to show that $\vec v \in S^1\mathcal C_x$. The proof involves a
discussion of holomorphic arcs that have $\vec v$ as a tangent line.
To this end, observe that it follows from the definition of
$\locus^2(x)$, and from Lemma~\ref{lem:tanofopen} that there exists a
holomorphic arc $\gamma : \Delta \to \locus^2(x)$ with $\gamma(0)=x$
such that
\begin{enumerate}
\item The tangent line to $\gamma$ at $0 \in \Delta$ is $\vec v$ and
  
\item \ilabel{i:gs} for every $s \in \Delta$, there exists a point
  $g_s \in \locus^1(x)$ such that $X$ contains a line through $g_s$
  and $\gamma(s)$, and a line through $g_s$ and $x$.
\end{enumerate}

By Theorem~\ref{thm:main1}, there is nothing to show if
$\gamma(\Delta) \subset \locus^1(x)$. We can thus assume that
$\gamma(\Delta) \not \subset \locus^1(x)$. In particular, we assume that
$g_s$ is not constantly equal to $x$.

The choice of $g_s$ in Item~\iref{i:gs} needs of course not be unique.
However, after a finite base change, replacing $\Delta$ by a finite
cover, we can assume without loss of generality that we can realize
$g_s$ as a function of $s$, i.e.~that we have an arc
$$
g: \Delta \to \locus^1(x).
$$

Although we may have $g(s) = x$ for a discrete set of $s \in \Delta$,
we still find a well-defined holomorphic arc
$$
\begin{array}{rccc}
G: & \Delta & \to & \mathcal K_x \\
& s & \mapsto & \text{line through $x$ and $g(s)$.}
\end{array}
$$
Likewise, we find a holomorphic arc
$$
\begin{array}{rccc}
F: & \Delta & \to & \mathcal K \\
& s & \mapsto & \text{line through $g(s)$ and $\gamma(s)$ }.
\end{array}
$$

\begin{notation}
  Let $\ell_G \in \mathcal K_x$ and $\ell_F \in \mathcal K$
  be the limiting lines $\ell_G = G(0)$, $\ell_F = F(0)$.
\end{notation}

\begin{rem}
  By definition, for all $s \in \Delta$, we have $\gamma(s) \in F(s)$.
  Since $\gamma(0) = x$, it follows that the limiting line $\ell_F$
  contains $x$, i.e.~$\ell_F \in \mathcal K_x$.
\end{rem}

\subsubsection*{Step 2: Limits of the Zariski tangent spaces of $F(s)
  \cup G(s)$}

For general $s \in \Delta$, the lines $F(s)$ and $G(s)$ do not
coincide, and the Zariski tangent space of the reducible conic $F(s)
\cup G(s)$ at the point of intersection is spanned by the direction
vectors of $G(s)$ and $F(s)$. As before, this yields a holomorphic map
$$
\begin{array}{rccc}
P: & \Delta & \to & \text{Grassmannian of 2-dimensional planes in $\P_N$} \\
& s & \mapsto & \text{plane that contains $G(s)$ and $F(s)$}
\end{array}
$$

\begin{rem}\label{rem:tP}
  Although the plane $P(s)$ is generally not contained in $X$, it
  intersects $X$ tangentially at $g(s)$, i.e.~$T_{P(s)}|_{g(s)}
  \subset T_X|_{g(s)}$. It follows directly from the definition that
  the line $\Lambda_s := \P(T_{P(s)}|_{g(s)}^\vee)$ is contained in
  $S^1 \mathcal C$ for all $s \in \Delta$.
\end{rem}

\begin{lem}\label{lem:tginP}
  The tangent line to the arc $\gamma$ at $0 \in \Delta$ is contained
  in the limiting plane $P(0)$.
\end{lem}
\begin{proof}
  The tangent line to the arc $\gamma$ at $0 \in \Delta$ is given by
  the zero-dimensional scheme $Q_0$ of length two that is the
  scheme-theoretic limit of the sequence of double points $Q_s := \{x,
  \gamma(s)\}$ as $s \to 0$. Since $Q_s \in P(s)$ for general $s \in
  \Delta$, the properness of the Hilbert-scheme implies that $Q_0 \in
  P(0)$.
\end{proof}

\subsubsection*{Step 3: End of proof if $\ell_G \not = \ell_F$}

Since $\gamma(0) = x$, the assumption $\ell_G \not = \ell_F$ implies
that $g(0) = \gamma(0) = x$ and $P(0)$ is spanned by the two lines
$\ell_G$ and $\ell_F$ through $x$. By Remark~\ref{rem:tP} and
Lemma~\ref{lem:tginP}, the tangent line to $\gamma$ at $0 \in \Delta$
is contained in $S^1 \mathcal C$. Remark~\ref{rem:closednessOfC} of
page~\pageref{rem:closednessOfC} implies that $S^1 \mathcal C|_x = S^1
\mathcal C_x$, and the proof of Theorem~\ref{thm:main3} is finished.

\subsubsection*{Step 4: End of proof if $\ell_G = \ell_F$}

Recall from Fact~\ref{fact:mostcurvesstd} that the line $\ell_G$ is
free.  We can therefore apply Proposition~\ref{prop:smx} which asserts
that $\mathcal C|_{g(0)}$ is smooth at $\P(T_{\ell_G}|_{g(0)}^\vee)$
and that the projective tangent space is $\mathbb P\left(\left.
    T_X|_{\ell_G}^+ \right|_{g(0)}^\vee\right)$.

We have seen in Remark~\ref{rem:tP} that $\Lambda_s =
\P(T_{P(s)}|_{g(s)}^\vee)$ is a secant line of $\mathcal C|_{g(s)}$
for all $s \in \Delta$. The assumption that $\ell_G = \ell_F$ implies
that the limiting line $\Lambda_0$ is tangent to $\mathcal C|_{g(0)}$.
The line $\Lambda_0$ is therefore tangent to $\mathcal C|_{g(0)}$ at
the point $\P(T_{\ell_G}|_{g(0)}^\vee)$.

As a next step, apply Lemma~\ref{lem:splitOfStd} to $\ell_G$ to see
that $\Lambda_0$ determines a positive subbundle $F$ of rank 2 in
$T_X|_{\ell_G}$.  Although the 2-dimensional plane $P(0)$ will
generally not be contained in $X$, it follows from the uniqueness
assertion of Lemma~\ref{lem:splitOfStd} that $P(0)$ is tangent to $X$
along the line $\ell_G$, and that $F$ is the restriction of the
tangent bundle $T_{P(0)}|_{\ell_G}$. The tangent space of $P(0)$ at
$x$ is therefore contained in the positive part $T_X|_{\ell_G}^+$, and
another application of Proposition~\ref{prop:smx} yields that the line
$\P(T_{P(0)}|_x^\vee) \subset \P(T_X|_x^\vee)$ is tangent to $\mathcal
C_x$.  Lemma~\ref{lem:tginP} then implies that the tangent line to the
arc $\gamma$ at $0 \in \Delta$ is contained in the tangent variety to
$\mathcal C_x$, in particular, that it lies in $S^1 \mathcal C_x$.
Theorem~\ref{thm:main3} is then shown. \qed

\section{Examples}
\label{chap:examples}

\subsection{Complete intersections}

Let $X \subset \P_N$ be a smooth complete intersection of
hypersurfaces of degrees $d_1, \ldots, d_m$. We have 
$$
-K_X = {\mathcal O}_X(p+2) \text{\quad where \quad} p = N - 1 -
\sum_{i=1}^m d_i.
$$
Let us assume that $p>0$. Then for each general $x \in X$, there
exists a line on $X$ through $x$ and the space ${\mathcal K}_x$ of
such lines is a complete intersection of dimension $p>0$. In fact, the
defining equations of the complete intersection ${\mathcal K}_x \cong
{\mathcal C}_x \subset \P(T_X|_x^\vee)$ are given by derivatives at
$x$ of the defining equations of $X$. In particular, $\mathcal C_x$ is
a nondegenerate, smooth complete intersection of dimension $p>0$.  It
follows that ${\mathcal K}_x$ is irreducible and there exists a unique
irreducible component ${\mathcal K}$ of the set of lines on $X$ whose
members cover $X$.

\begin{prop}\label{prop:4}
  Let $Y \subset \P_m$ be a non-degenerate smooth complete
  intersection of dimension $p >0$.  Then $Y$ has no secant defect,
  i.e., $\dim S^kY = \min\{m, (k+1)(p+1)-1 \}$ for all $k$.  
\end{prop}

\begin{proof}
  Otherwise there exists a hyperplane in $\P_m$ tangent to $Y$ along a
  positive-dimensional subvariety by \cite[V.1.5]{Zak93}. This is
  impossible by \cite[7.5]{FL81}.
\end{proof}

\begin{prop}\label{prop:5}
  Let $X \subset \P_N$ be an $n$-dimensional smooth complete
  intersection with $-K_X = {\mathcal O}_X(p+2)$ and $p>0.$ Then for
  each $k>0$, the dimension $d_k$ equals
  $$
  d_k= \min\{ n-1, \, pk + k-1 \} +1.
  $$
\end{prop}

\begin{proof}
  Let ${\mathcal C}_x \subset \P(T_X|_x^\vee)$ be the variety of
  tangents to lines through $x$ for a general $x \in X$, and suppose
  $pk+k-1 \leq n-1$. Then by Proposition~\ref{prop:1},
  Theorem~\ref{thm:main2} and Proposition~\ref{prop:4},
  $$
  pk+k-1 = \dim S^{k-1} {\mathcal C}_x \leq \dim {\mathcal C}_x^k
  \leq k (p+1)-1.
  $$
  Thus 
  $$
  d_k-1= \dim {\mathcal C}_x^k = \dim S^{k-1} {\mathcal C}_x =
  pk+k-1.
  $$
  If $pk+k-1 \geq n-1$, the equality $S^k {\mathcal C}_x =
  {\mathcal C}_x^k = \P(T_X|_x^\vee)$ is obvious.  
\end{proof}

\subsection{Hermitian symmetric spaces}

Let $X= G/P$ be an irreducible Hermitian symmetric space. The action
of the isotropy subgroup $P$ at a base point $x \in X$ on the tangent
space $\P(T_X|_x^\vee)$ is irreducible. This action has exactly $r$
orbits where $r$ is the rank of the Hermitian symmetric space. There
exists a minimal equivariant embedding $X \subset \P_N$ such that the
set of lines of $\P_N$ lying on $X$ forms a family ${\mathcal K}$ that
satisfies the Assumption~\ref{ass:dagger}.  The highest weight orbit
of the isotropy action on $\P(T_X|_x^\vee)$ is then precisely the
variety of minimal rational tangents ${\mathcal C}_x$ for ${\mathcal
  K}$. The possibilities are listed in Table~\ref{tab:typesofhermss}.

\begin{table}[tbp]
  \small

  \begin{tabular}{ccccccc} \hline
    Type & $G$ & $X$ & $\dim(X)$ & rank & $\dim({\mathcal C}_x)$ &  ${\mathcal C}_x$ \\
    \hline I & $SL(a+b)$ & $G(a,b)$ & $ab$ & $\min(a,b)$ & $a+b-2$ &
    $\P_{a-1} \times \P_{b-1}$ \\
    II & $SO(2m)$ & $QG(m)$ & $\frac{m(m-1)}{2}$& $[\frac{m}{2}]$ &
    $2m-4$ & $G(2, m-2)$ \\ III & $Sp(m)$ & $LG(m)$ &
    $\frac{m(m+1)}{2}$ & $m$ & $m-1$ & $v_2(\P_{m-1})$ \\
    IV & $SO(m+2)$ & ${\mathbb Q}_m$ & $m$ & 2 & $m-2$ & ${\mathbb Q}_{m-2}$ \\
    V & $E_6$ & ${\mathbb O}\P_2
    \otimes_{\mathbb R} {\mathbb C}$ & 16 & 2& 10 & $QG(5)$  \\
    VI & $E_7$ & no classical name & 27 & 3 & 16 & ${\mathbb O}\P_2
    \otimes_{\mathbb R} {\mathbb C}$  \\ \hline \ \\
  \end{tabular}

  \begin{tabular}{lcl}
    $G(a,b)$ & $\ldots$ & Grassmannian of $a$-dim.~subspaces in
    an $(a+b)$-dim.~vector space \\
    $QG(m)$ & $\ldots$ & quadric Grassmannian of $m$-dimensional isotropic
    subspaces in a \\
    & & $2m$-dimensional orthogonal vector space\\
    $LG(m)$ & $\ldots$ & Lagrangian Grassmannian of a $2m$-dimensional
    symplectic vector space \\
    $v_2(\P_{m-1})$ & $\ldots$ & second Veronese embedding \\
    ${\mathbb O}\P_2$ & $\ldots$ & the octahedral plane \\
    $[\cdots]$ & $\ldots$ & round-down.
  \end{tabular}
  
  \medskip

  \caption{Types of Hermitian Symmetric Spaces}
  \label{tab:typesofhermss}
\end{table}

Let us describe the higher secant varieties of ${\mathcal C}_x$, case
by case. The descriptions are fairly standard, see e.g.
\cite[1.3.6]{Flenner-OCarrol-Vogel}.

\subsubsection*{Type I}

${\mathcal C}_x \subset \P(T_X|_x^\vee)$ is isomorphic to the Segre
variety $\P_{a-1} \times \P_{b-1} \subset \P_{ab-1}$.  Its affine cone
corresponds to the set of elements of rank $\leq 1$ in the space of
$(a \times b)$-matrices. Its $k$-secant variety is the set of elements
of rank $\leq k+1$. Since there are $\min(a,b)$ orbits in
$\P(T_X|_x^\vee)$, each orbit closure is one of the higher secant
varieties of ${\mathcal C}_x$. Their dimensions are
$$
\dim S^k {\mathcal C}_x = (a+b)(k+1) -(k+1)^2 -1 \mbox{ for } 1 \leq k \leq
\min(a,b)-1.
$$

\subsubsection*{Type II}

${\mathcal C}_x \subset \P(T_X|_x^\vee)$ is the Plücker embedding of
$G(2, m-2)$. Its affine cone corresponds to the set of elements of
rank $\leq 2$ in the space of alternating $(m\times m)$-matrices. Its
$k$-secant variety is the set of elements of rank $\leq 2(k+1)$. Since
there are $[\frac{m}{2}]$ orbits in $\P(T_X|_x^\vee)$, each orbit
closure is one of the higher secant varieties of ${\mathcal C}_x$.
Their dimensions are
$$
\dim S^k {\mathcal C}_x = \frac{m(m-1)}{2} -
\frac{(m-2k-2)(m-2k-3)}{2} -1 \mbox{ for } 1 \leq k \leq
\left[\frac{m}{2}\right]-1 .
$$

\subsubsection*{Type III}

${\mathcal C}_x \subset \P(T_X|_x^\vee)$ is the second Veronese
embedding of $\P_{m-1}$. Its affine cone corresponds to the set of
elements of rank $\leq 1$ in the space of symmetric $(m \times
m)$-matrices. Its $k$-secant variety is the set of elements of rank
$\leq k+1$. Since there are $m$ orbits in $\P(T_X|_x^\vee)$, each orbit
closure is one of the higher secant varieties of ${\mathcal C}_x$.
Their dimensions are
$$
\dim {S^k\mathcal C}_x = \frac{m(m+1)}{2} - \frac{(m-k)(m-k-1)}{2}
-1 \mbox{ for } 1 \leq k \leq m-1.
$$

\subsubsection*{Type IV}

There are only two orbits in $\P(T_X|_x^\vee)$, the hypersurface
${\mathcal C}_x$ and its complement. It follows that $S^1{\mathcal
  C}_x = \P(T_X|_x^\vee)$.

\subsubsection*{Type V}

There are only two orbits in $\P(T_X|_x^\vee)$, ${\mathcal C}_x$ and
its complement. It follows that $S^1 {\mathcal C}_x =
\P(T_X|_x^\vee)$.

\subsubsection*{Type VI}

${\mathcal C}_x \subset \P(T_X|_x^\vee)$ is the $E_6$-Severi variety.
By \cite[p. 59, E]{Zak93}, $S^1 {\mathcal C}_x$ is then a cubic
hypersurface. There are only three orbit closures in
$\P(T_X|_x^\vee)$, namely ${\mathcal C}_x$, the cubic hypersurface
${S^1\mathcal C}_x$, and $S^2 {\mathcal C}_x = \P(T_X|_x^\vee)$.

\subsubsection*{The length of Hermitian symmetric spaces}
In summary, we have obtained the following which must be well-known to
experts.

\begin{prop}\label{prop:6}
  Let $X$ be an irreducible Hermitian symmetric space of rank $r$.
  Then there are exactly $r$ projective subvarieties in
  $\P(T_X|_x^\vee)$ invariant under the isotropy action, namely
  ${\mathcal C}_x$ and its higher secant varieties.
\end{prop}

The spanning dimension of chains of minimal rational curves can be
computed now.

\begin{lem}\label{lem:8}
  The length of an irreducible Hermitian symmetric space of rank $r$
  is at least $r$.  
\end{lem}

\begin{proof}
  This follows from Proposition~\ref{prop:3} and the well-known fact
  that there exists a section of the ample generator of $\Pic(X)$
  which vanishes to order $r$ at some point. In fact, this section is
  the compactifying divisor of the Harish-Chandra embedding of
  $\mathbb C^n$ in $X$.
\end{proof}

\begin{prop}\label{prop:7}
  Let $X$ be an irreducible Hermitian symmetric space. Then for each
  $k$, $S^k{\mathcal C}_x = {\mathcal C}_x^{k+1}$.  
\end{prop}

\begin{proof}
  Since $S^k {\mathcal C}_x \subset {\mathcal C}^{k+1}_x$, and both
  $S^k {\mathcal C}_x$ and ${\mathcal C}^{k+1}_x$ are invariant
  subvarieties of the isotropy action, Proposition~\ref{prop:6}
  implies that if there exists a number $j$ with $S^j {\mathcal C}_x
  \ne {\mathcal C}^{j+1}_x$, then ${\mathcal C}^{r-1}_x =
  \P(T_X|_x^\vee)$. This is a contradiction to the preceding
  Lemma~\ref{lem:8}.
\end{proof}

As a direct consequence of Propositions~\ref{prop:6} and \ref{prop:7},
we obtain the following.

\begin{cor}
  The length of an irreducible Hermitian symmetric space of rank $r$
  is precisely $r$.  
\end{cor}

\subsection{Homogeneous contact manifolds}
\label{sect:contactmfs}

Let $X=G/P$ be a homogeneous contact manifold of dimension $2m+1$ with
$b_2(X)=1$. The action of the isotropy group $P$ at a base point $x
\in X$ on the tangent space $T_X|_x$ has a unique irreducible subspace
$D_x$ of dimension $2m$.  It is known that $X$ can be embedded into a
projective space by the ample generator ${\mathcal O}_X(1)$ of
$\Pic(X)$. Assume that $X$ is different from $\P_{2m+1}$. Then $-K_X =
{\mathcal O}_X(m)$. Let $\mathcal K$ be the set of lines lying on $X$.
Then ${\mathcal C}_x \subset \P(D|_x^\vee)$ is the highest weight
orbit of the isotropy action of $P$ on $\P(D|_x^\vee)$.
Table~\ref{tab:2}, which is taken from \cite[p.~353]{Hwa00},
summarizes the possibilities.
\begin{table}[tbp]
  \small
  \begin{tabular}{cccc}
    \hline $G$ & $m $ & ${\mathcal C}_x$ & embedding ${\mathcal C}_x
    \subset \P(D|_x^\vee)$ \\
    \hline $SO(m+4)$ & $m \geq 3$& $\P_1 \times {\mathbb Q}_{m-2}$ & Segre\\
    $G_2$ & 2 & $\P_1$ & third Veronese\\ 
    $F_4$ & 7 & $LG(3)$ & m.e. \\
    $E_6$ & 10 & $G(3,3)$ & m.e. \\
    $E_7$ & 16 & $QG(6)$ & m.e. \\
    $E_8$ & 28& Type VI & m.e.\\
    \hline \ \\
  \end{tabular}
  
  \begin{tabular}{lcl}
    m.e. & \ldots & minimal equivariant embedding of irreducible Hermitian
    symmetric spaces
  \end{tabular}

  For other notation, see the cutline of Table~\ref{tab:typesofhermss}
  
  \medskip
  
  \caption{Homogeneous contact manifolds}
  \label{tab:2}
\end{table}
In all cases, 
\begin{align*}
  \dim {\mathcal C}_x & = m-1 \\
  \dim \P(D|_x^\vee) & = 2m-1 = 2 \cdot \dim {\mathcal C}_x +1
\end{align*}
By \cite[III.2.1, III.1.4, and III.1.7]{Zak93}, we have
$$
S^1 {\mathcal C}_x = \P(D|_x^\vee) = {\mathcal C}^2_x.
$$
In particular, $\locus^2(x)$ is a hypersurface in $X$ and $\locus^3(x)
= X$ by Proposition~\ref{prop:nadel}. As a consequence we obtain that
$$
S^2{\mathcal C}_x = \P(D|_x^\vee) \neq \P(T_X|_x^\vee) = {\mathcal
  C}^3_x,
$$
which shows that the inclusion in Theorem~\ref{thm:main2} can
sometimes be strict.

\subsection{Moduli space of stable bundles of rank 2 over a
curve}

Throughout the present section, let $C$ be a smooth projective curve
of genus $g \geq 4$ and $X := {\mathcal SU}_C(2, i)$ be the moduli
space of stable vector bundles of rank 2 with a fixed determinant of
degree $i \in \{0,1\}$. 

\begin{fact}[{\cite[sect.~1.4.8]{Hwa00b}}]\label{fact:mod}
  The space $X$ is a smooth quasi-projective variety of dimension
  $3g-3$ which is projective if $i=1$. There exists a very ample line
  bundle $L$ which generates the Picard group, $\Pic(X) \cong {\mathbb
    Z} \cdot L$. Moreover, we have $-K_X = 4 \cdot L$ when $i=0$ and
  $-K_X = 2 \cdot L$ when $i=1$.  The moduli space $X$ is covered by a
  family $\mathcal K$ of compact rational curves, called \emph{Hecke
    curves}, that satisfy Assumption~\ref{ass:dagger}. These curves
  intersect $-K_X$ with multiplicity $4$.
\end{fact}

\begin{rem}
  In the previous sections we assumed that $X$ is a smooth projective
  variety. We remark that Definitions~\ref{defn:VMRT} and \ref{def:Cx}
  work without change for the quasi-projective manifold $SU_C(2,0)$.
\end{rem}

\begin{fact}[{\cite[prop.~11]{Hwa00b}}]\label{fact:Hecke}
  Let $x \in X$ be a general point, and $W$ be the corresponding
  vector bundle on $C$. Let $Y := \P(W)$ be the associated ruled
  surface over $C$, and $\pi: Y \to C$ be the canonical projection.
  Then the linear system of the line bundle $H := 2\pi^*(K_C) - K_Y$
  is basepoint-free, and determines a finite, generically injective
  map
  $$
  \phi_H : Y \to \P_{3g-4}
  $$
  Further, there exists an identification $\P(T_X|_x^\vee) \cong
  \P_{3g-4}$ such that the variety of minimal rational tangents
  $\mathcal C_x \subset \P(T_X|_x^\vee)$ is exactly identified with
  the image of $\phi_H$. In particular,
  \begin{enumerate}
  \item $\dim \mathcal C_x = 2$,
  \item $\mathcal C_x$ is not rational and not isomorphic to a cone
    over a curve, and
  \item $\mathcal C_x$ is not linearly degenerate.
  \end{enumerate}
\end{fact}

\subsubsection{The variety of minimal rational tangents}

In this section, we use Fact~\ref{fact:Hecke} and classic results of
projective geometry to show that the variety of minimal rational
tangents on $X$ has no secant defect.

\begin{prop}\label{prop:9}
  Let $x \in X$ be a general point, and $W$ be the associated stable
  bundle of rank 2 over the curve $C$ of genus $g\geq 4$. Let
  ${\mathcal C}_x$ be the variety of tangents with respect to the
  space $\mathcal K$ of Hecke curves.  Then ${\mathcal C}_x$ has no
  secant defect, i.e.,
  $$
  \dim S^k {\mathcal C}_x = \min \{ 3g-4, \, 3(k+1)-1\}.
  $$
  In particular, $d_k= \min\{3g-3, 3k \}$ for each $k$.  
\end{prop}

The proof of Proposition~\ref{prop:9} is based on the following lemma.
Its proof is implicitly contained in \cite[proof of prop.~12]{Hwa00b}.
However, since the notation there is somewhat different, we reproduce
the full proof for the readers' convenience.

\begin{lem}\label{lem:610}
  In the setup of Proposition~\ref{prop:9}, suppose there exists a
  linear subspace $\P_s \subset \P(T_X|_x^\vee)$ such that the tangent
  plane to ${\mathcal C}_x$ at a general point has non-empty
  intersection with $\P_s$. Then $s \geq g-2$.
\end{lem}
\begin{proof}
  Choose a general complementary linear subspace $\P_{3g-5-s} \subset
  \P(T_X|_x^\vee)$ and consider the projection $\psi$ from $\P_s$ to
  $\P_{3g-5-s}$. We obtain a diagram
  $$
  \xymatrix{ {Y} \ar[rr]^(.4){\phi_H}
    \ar[d]_{\txt{\scriptsize $\pi$ \\ \scriptsize $\P_1$-bundle}} & &
    {\P(T_X|_x^\vee)}
    \ar@{-->}[d]^{\txt{\scriptsize $\psi$ \\ \scriptsize projection
        from $\P_s$}} \\
    {C} & & {\P_{3g-5-s}} }
  $$
  where the ruled surface $Y$ and the morphism $\phi_H$ are the
  same as in Fact~\ref{fact:Hecke} above.
  
  \subsubsection*{The tangent map of $\psi$}
  If $y \in \mathcal C_x$ is a general point, and $T$ its projective
  tangent space, the fact that ${\mathcal C}_x \not \subset \P_s$
  implies that the intersection of a general tangent plane to
  ${\mathcal C}_x$ and $\P_s$ can only have dimension 0 or 1. In these
  cases, the rank of the tangent map $T(\psi|_{\mathcal C_x})$ at $y$
  is 1 or 0, respectively. Since $y \in \mathcal C_x$ is a general
  point, this means that the image of $\mathcal C_x$ under the
  projection is either a curve $Z$ or a point in $\P_{3g-5-s}$.
  
  \subsubsection*{Reduction to the case of 1-dimensional fibers}
  If the intersection has dimension 1, then the projection $\psi$
  sends the tangent space at a general point of ${\mathcal C}_x$ to a
  single point. Thus the projection sends ${\mathcal C}_x$ to a point.
  This implies that ${\mathcal C}_x$ is contained in a fiber of the
  projection map, i.e.~in linear subspace $\P_{s+1} \supset \P_s$,
  which contradicts the nondegeneracy of ${\mathcal C}_x$ unless $s+1
  \geq 3g-4$ in which case the proof is finished.
  
  We may therefore assume that the projection $\psi$ maps ${\mathcal
    C}_x$ to a curve.

  \subsubsection*{Fibers of $\psi\circ \phi_H$}
  Suppose that the $\phi_H$-image of a generic $\pi$-fiber $f
  \subset Y$ dominates $Z$. Since 
  $$
  \deg Z \leq \deg \phi_H(f) \leq H.f = 2,
  $$
  $Z$ is then either a line or a conic, and is contained in a
  2-dimensional linear plane. This would imply that ${\mathcal C}_x$
  is contained in a linear subspace $\P_{s+3} \supset \P_s$ which
  contradicts the nondegeneracy of ${\mathcal C}_x$ unless $s+3 \geq
  3g-4$ in which case the proof is finished.
  
  We may therefore assume that $\psi$ contracts general fibers of
  $\pi$ to a point. In particular, we assume that $\phi_H(f)$ is
  contained in a linear space $\P_{s+1}$ that contains $\P_s$ as a
  hyperplane, and intersects $\P_s$ twice, counting multiplicity. In
  other words, the cycle-theoretic preimage $\phi_H^*(\P_s)$ contains
  a double section $D$.
  
  \subsubsection*{Discussion of linear systems on $Y$}
  Let $\Xi \subset |H|$ be the subsystem of dimension $3g-5-s$ that
  defines the map $\psi\circ \phi_H$. Since $D \subset Y$ is a fixed
  component of $\Xi$, we have
  \begin{equation}
    \label{eq:linser}
    3g-5-s = \dim |\Xi| \leq \dim |H-D|    
  \end{equation}
  Using the notation of \cite[V.2.8.1]{Ha77}, there exists an integer
  $a$ and line bundles $\mathfrak a, \mathfrak b \in \Pic(C)$ such
  that
  \begin{align*}
    D & \equiv 2 C_0 + \mathfrak a f & \deg \mathfrak a & = a \\
    H & \equiv 2 C_0 + \mathfrak b f & \deg \mathfrak b & = 2g-2+e.
  \end{align*}
  Recall that since the bundle $W$ is stable, we have $e < 0$. In
  particular, by \cite[V.~prop.~2.21]{Ha77}, the minimal section $C_0$
  is ample. As a consequence we have the following inequalities
  \begin{align*}
    -C_0 \cdot D = 2e - a & < 0 && \text{since $C_0$ is ample, and} \\
    -g & \leq e && \text{by \cite{Nagata70}.} 
  \end{align*}
  Together they imply that
  $$
  \deg (\mathfrak b - \mathfrak a) = 2g-2+e-a \leq 3g-3,
  $$
  and an application of Riemann-Roch and of Clifford's theorem,
  \cite[p.~343]{Ha77}, yields
  $$
  \dim |H-D| = \dim |\mathfrak b - \mathfrak a| \leq
  \max\left(\frac{3g-3}{2},\, 2g-3\right)= 2g-3.
  $$
  Combined with the inequality~\eqref{eq:linser}, we get $g-2 \leq
  s$ as desired. This ends the proof of Lemma~\ref{lem:610}.
\end{proof}

With the aid of Lemma~\ref{lem:610}, Proposition~\ref{prop:9} follows
easily from the classification \cite{CC02} of weakly defective
projective surfaces. We remark that there is a misprint in this paper
which affects us here. The last sentence of \cite[thm.~1.3(ii)]{CC02}
should read ``The minimal such $s$ is characterized by the property
that $X$ is $(s+1)$-defective but not $s$-defective, [\ldots]''

\begin{proof}[Proof of Proposition~\ref{prop:9}]
  We argue by contradiction and assume that ${\mathcal C}_x$ has a
  positive secant defect. Let $k \geq 1$ be the minimal number such
  that ${\mathcal C}_x$ is $k$-defective. It follows from
  Fact~\ref{fact:Hecke} and from \cite[thm.~1.3 (i) and (ii)]{CC02},
  that
  \begin{equation}
    \label{eq:wdef}
     k \leq g-2,    
  \end{equation}
  and that ${\mathcal C}_x$ is contained in a $(k+1)$-dimensional cone
  over a curve whose vertex is a linear space $V \subset
  \P(T_X|_x^\vee)$ of dimension $\dim V = k-1$. By
  \cite[prop.~4.1]{CC02}, this description implies the following: if
  $y \in {\mathcal C}_x$ is a general point, and $T_{{\mathcal C}_x,y}
  \subset \P(T_X|_x^\vee)$ the projective tangent space to ${\mathcal
    C}_x$ at $y$, then $V \cap T_{{\mathcal C}_x,y} \ne \emptyset$.
  In this setup, an application of Lemma~\ref{lem:610} shows that
  \begin{equation}
    \label{eq:dimv}
    k-1 = \dim V \geq g-2
  \end{equation}
  The inequalities~\eqref{eq:wdef} and~\eqref{eq:dimv} are clearly
  contradictory. We have thus seen that $\mathcal C_x$ has no secant
  defect.
  
  Theorem~\ref{thm:main2} and Proposition~\ref{prop:1} then yield that
  $d_k = \min( 3g -3,3k)$.
\end{proof}

As a direct corollary, we have the following bound on the
multiplicities of divisors on ${\mathcal SU}_C(2, i)$.

\begin{cor}\label{cor:c11}
  Let $x\in {\mathcal SU}_C(2, i)$ be a general point, let $L \in
  \Pic({\mathcal SU}_C(2, i))$ be the ample generator of the Picard
  group, and $D \in |mL|$, $m \geq 1$ any divisor. Then the
  multiplicity of $D$ at $x$ is bounded as follows
  \begin{equation*}
    \mult_x (D) \leq \left\{ 
      \begin{array}{ll}
        m(g-1) & \text{if $\, i=0$} \\
        2m(g-1) & \text{if $\, i=1$}
      \end{array}
    \right.
  \end{equation*}
\end{cor}
\begin{proof}
  Let $\ell \subset X$ be a Hecke curve.
  
  If $i = 1$, then Proposition~\ref{prop:9} implies that the length of
  $X$ is $g-1$. Fact~\ref{fact:mod} implies that $\ell \cdot L = 2$,
  i.e.~$\ell \cdot D = 2m$, and the claim follows immediately from
  Proposition~\ref{prop:3}.
  
  If $i = 0$, we consider a smooth compactification $X \subset
  \overline X$, the associated compactification $D \subset \overline
  D$, and let $\overline{\mathcal K} \subset \RatCurves^n(\overline
  X)$ be the closure of $\mathcal K$ under the natural inclusion
  morphism $\RatCurves^n(X) \subset \RatCurves^n(\overline X)$. It is
  clear that the Assumptions~\ref{ass:dagger} hold for the family
  $\overline {\mathcal K}$ of curves on $\overline X$, and that the
  length of $\overline X$ with respect to $\overline{\mathcal K}$ is
  then also $g-1$. Since $\ell \cdot \overline D = \ell \cdot D =1$,
  the claim follows again from Proposition~\ref{prop:3}.
\end{proof}

\subsubsection{Acknowledgments}
We would like to thank C.~Ciliberto and F.~L.~Zak for kindly answering
questions regarding higher secants and M.~Popa for a discussion on
Corollary~\ref{cor:c11}. We are grateful to the referee for careful
reading and for suggestions.

 \end{document}